\documentclass[12pt]{article}
\newtheorem{thm}{Theorem}
\newtheorem{prop}[thm]{Proposition}
\newtheorem{cor}[thm]{Corollary}
\newtheorem{lemma}[thm]{Lemma}

\def\proof{{\sc Proof. }}

\usepackage{amsfonts}
\usepackage{amssymb}

\def\CP{\mathbb{C}{\rm P}}
\def\Re{{\rm Re}}

\def\C{\mathbb{C}}
\def\H{\mathbb{H}} 
\def\R{\mathbb{R}}
\def\Fc{{\cal F}}
\def\Gc{{\cal G}}
                  
\def\Box{\square}

\def\geq{\geqslant}
\def\<{\langle}
\def\>{\rangle}

\title{K\"ahler metrics whose geodesics are circles}
\author{V. Timorin
\thanks{Partially supported by CRDF RM1-2086}}
\date{}

\begin{document}
\maketitle

{\small We classify all K\"ahler metrics in an open subset of $\C^2$ whose
real geodesics are circles. All such metrics are equivalent (via complex
projective transformations) to Fubini metrics (i.e. to Fubini-Study metric
on $\CP^2$ restricted to an affine chart, to the complex hyperbolic
metric in the unit ball model or to the Euclidean metric).}

\section*{Introduction}

All Riemannian metrics in an open subset of $\R^2$ or $\R^3$ whose
geodesics are arcs of circles are classical, i.e., isometric
to Euclidean, Riemann or Lobachevsky geometries. This was proved
by A. Khovanskii \cite{Kh} in dimension 2 and by F. Izadi \cite{Iz} in
dimension 3. But in dimension 4 this is wrong. There are remarkable
K\"ahler metrics whose real geodesics are circles --- Fubini metrics
(see Appendices 1 and 2).

Our main result is as follows:

\begin{thm}
Consider a K\"ahler metric in an open subset of $\C^2$ such that all
geodesics are parts of circles (or straight lines). Then this metric is
(up to a complex projective transformation) some Fubini metric. 
\end{thm}

In the next Section we will prove this result. Then we mention
(without proof) a local geometric classification of complete families of
circles that are point-wise rectifiable by means of complex projective
transformations. 

{\em Acknowledgements.}
I am grateful to A.G. Khovanskii and R.W. Sharpe for useful discussions. 

\section*{Proof of the main result}

For a definition of K\"ahler metrics see Appendix 2. 
Let $g$ be a K\"ahler metric in an open region $\Omega\subset\C^2$ such that
all geodesics with respect to $g$ are parts of circles. First note that
for any point $p\in\Omega$ the set of geodesics passing through $p$ coincides
with the image under the exponential map of the set of lines passing through
$p$.

The following theorem \cite{Tim} holds:

\begin{thm}
Fix some identification of $\R^4$ with the algebra $\H$ of quaternions.
Suppose that a local diffeomorphism $\Phi:(\R^4,0)\to(\R^4,0)$ takes
sufficiently many lines (through 0) in general position to circles.
If $d_0\Phi=id$, then the second derivative of $\Phi$ has the form
$x\mapsto A(x)x$ or $x\mapsto xA(x)$ where $A$ is some $\R$-linear map, and
the multiplication is in the sense of quaternions. 
\end{thm}

In particular, the bundle of geodesics at a point $p\in\Omega$ is given
by $p+xt+\frac 12 A(x)xt^2$ or $p+xt+\frac 12 xA(x)t^2$
where $t$ is a parameter, $x$ is the velocity vector and $A$ is some linear
map.
To fix the idea assume that the multiplication by $A(x)$ is from the left.
Now recall the following well-known fact (for a proof, see Appendix 2)

\begin{prop}
\label{exp-hol}
Exponential maps with respect to a K\"ahler metric are holomorphic
up to third order terms (i.e., their 2-jets are holomorphic).
\end{prop}

Therefore, the map $x\mapsto A(x)x$ must be holomorphic. Now we need

\begin{lemma}
Suppose that the map $x\mapsto A(x)x$ (where $A$ is some linear operator)
is holomorphic. Then $A(x)$ is complex linear and takes complex values only.
In other words, $A$ is a complex linear functional. 
\end{lemma}

\proof
In general, $A(x)=a(x)+b(x)i+c(x)j+d(x)k$ where $a,b,c,d$ are some linear
functionals on $\R^4$. Being holomorphic, the quadratic map $x\mapsto A(x)x$
must satisfy the condition $A(ix)(ix)=-A(x)x$.
From this condition it follows that
$$a(x)=b(ix),\quad b(x)=-a(ix),\quad c(x)=-d(ix),\quad d(x)=c(ix).$$
This means that the functionals $\alpha=a+bi$ and $\beta=d+ci$ must be 
holomorphic (i.e. complex linear). Since the map
$x\mapsto A(x)x=(\alpha(x)+k\beta(x))x$ is holomorphic, the map
$x\mapsto k\beta(x)x$ must be also holomorphic. If
we multiply $x$ by $\sqrt{i}$, then $k\beta(x)x$ gets multiplied by $-i$,
but by bilinearity it must be multiplied by $i$. Hence $\beta=0$ and
$A(x)=\alpha(x)$ is a complex linear functional. $\Box$

From this lemma it follows in particular that all geodesics of $g$
lie in complex lines (since velocities and accelerations are proportional
with some complex coefficient). Now we can use a complexified version
of Beltrami's theorem (it follows from the results of Bochner \cite{Bo},
Otsuki and Tashiro \cite{OT}; for a detailed discussion and a sketch of
a proof see Appendix 3):

\begin{prop}
\label{B}
Suppose that all germs of complex lines are totally geodesic surfaces
with respect to some Hermitian metric on a part of $\C^n$. Then this metric
is equivalent (up to a complex projective transformation) to a Fubini metric. 
\end{prop}

This concludes the proof of the main theorem.

\section*{Complex families of circles}

Consider a K\"ahler metric $g$ defined in an open region
$\Omega\subseteq\C^2$ and assume that all geodesics of $g$ are arcs of
circles.
We saw that the exponential map of $g$ at any point $p\in\Omega$ has the
form $x\mapsto p+x+\frac 12 A(x)x$ or $x\mapsto p+x+\frac 12 xA(x)$ up to
third order terms. Here $A$ is some complex
linear functional. Note that a complex projective transformation
$x\mapsto p+(1-\frac 12 A(x))^{-1}x$ or 
$x\mapsto p+x(1-\frac 12A(x))^{-1}$ has the same 2-jet and clearly takes all
lines to circles. Therefore, the images of lines (through 0) under the
above complex projective map are geodesics of $g$ (through $p$). This is
because a circle is determined by its velocity and acceleration at some
point and hence by the 2-jet of some rectifying diffeomorphism. Thus
the geodesics of $g$ are point-wise rectifiable by means of complex
projective transformations.

Suppose that we are given a family $\Fc$ of curves in $\C^2$ (by a curve
we mean a 1 dimensional closed submanifold).
Let us say that the family $\Fc$ is a {\em complete family of curves}
in an open subset $\Omega$ of $\C^2$ if through each point of $\Omega$
in each direction there goes a curve from $\Fc$. A family $\Fc$ is said
to be {\em rectifiable} at some point $p\in\Omega$ if there exists a
germ of diffeomorphism at $p$ that takes each curve from $\Fc$ passing
through $p$ to a straight line. A complete family $\Fc$ of curves is
called a {\em complex family of curves} in $\Omega$ if it is point-wise
rectifiable in $\Omega$ by means of local diffeomorphisms holomorphic up
to third order terms.
Complete complex families of curves generalize the notion of geodesics
with respect to a K\"ahler metric. As we saw a complex family of circles
is point-wise rectifiable by means of complex projective transformations.

We have a local classification of all complete complex families
of circles in $\Omega$. Up to a complex projective transformation
these are the following:

\begin{itemize}
\item Geodesics of Fubini metrics.
\item A family outside the unit ball.
\item Suspensions. 
\end{itemize}

Let us describe the last 2 examples in detail.

{\bf A family of circles outside the unit ball.} Inside the unit ball
we have a model for the complex hyperbolic plane. The metric is given in
coordinates $(z_1,z_2)$ by
$$ds^2=\frac{(dz_1d\bar z_1+dz_2d\bar z_2)(1-z_1\bar z_1-z_2\bar z_2)
+(dz_1\bar z_1+dz_2\bar z_2)(d\bar z_1z_1+d\bar z_2z_2)}
{(1-z_1\bar z_1-z_2\bar z_2)^2}$$ 
All geodesics are circles (see Appendix 2). Note that this metric makes
sense in the exterior of the unit ball as well. But in the exterior it will
be no more positive definite. Nevertheless, geodesics make sense and they
are circles. We get a complex family of circles that are not geodesics
with respect to a (positive definite) K\"ahler metric (this is not
obvious but not very hard to prove).

{\bf Suspensions.} Let $U$ be a domain in $\C$ and $\Fc$ a point-wise
rectifiable family of circles in $U$ (e.g. the set of geodesics in the
Poincar\'e half-plane). Note that any point-wise rectifiable family in
dimension 2 is complex (this follows from the result of Khovanskii
\cite{Kh}). We are going to define a complex family $\Gc$ of
circles in $U\times\C$ that will be called the {\em suspension} of $\Fc$.
Each circle from $\Gc$ must lie in some complex line. We will define 
$\Gc$ on every complex line separately, and then prove that $\Gc$ is
point-wise rectifiable. 
    
Take any complex line $L$ in $\C^2$. Then the projection $\pi$ from
$L\cap (U\times\C)$ to $U$ either maps everything to a point or is a linear
conformal one-to-one map. In the first case (i.e. when $L$ is ``vertical'')
define $\Gc$ on $L$ as the set of all real lines in $L$. In the second case
the map $\pi^{-1}$ clearly takes circles to circles. So define $\Gc$ on $L$
as the preimage of the set of all circles in $U$ under the projection $\pi$.

Let us prove that the family $\Gc$ thus constructed is a complex family
of circles.
Take a point $a\in U\times\C$. Suppose that $\Fc$ can be rectified at the
point
$\pi(a)$ by some complex projective map $P=L_1/L_2$ where $L_1$ and $L_2$ are
affine functions. It is easy to see that the map $(z,w)\mapsto (P(z),w/L_2(z))$ 
rectifies the family $\Gc$ at $a$.

The proof of the above classification is not very hard. Nevertheless we
will not give it here. 

\section*{Appendix 1: Hermitian and K\"ahler metrics}

Consider a Riemannian metric $g$ in an open subset $\Omega$ of $\C^n$.
This metric is called {\em Hermitian} if it is stable under the multiplication
by $i$, i.e. $g(ix,iy)=g(x,y)$ for any 2 vectors $x$ and $y$ at the
same point. With a Hermitian metric $g$ one associates the differential
$(1,1)$-form $\omega(x,y)=g(ix,y)$ and a sesquilinear form
(Hermitian inner product) $\<X,Y\>=g(X,Y)-i\omega(X,Y)$.
A metric $g$ is said to be {\em K\"ahler} if $d\omega=0$. 

Let $\nabla^0$ be the standard (flat) connection on $\C^n$.
Denote by $\nabla$ the Levi-Civita connection of $g$.
Then for each pair of vector fields $X$ and $Y$ on $\Omega$ we have
$\nabla_X Y=\nabla^0_X Y+\Gamma(X,Y)$ where $\Gamma$ is a symmetric
$\R$-bilinear form at each point. The form $\Gamma$ is called {\em
Christoffel form}. In particular, the value of $\Gamma(X,Y)$
at a point $p$ depends only on the values of $X$ and $Y$ at $p$
(not on their derivatives). Let us recall the following fact:

\begin{prop}
\label{K-cri}
A metric $g$ is K\"ahler if and only if the corresponding covariant
differentiation 
is complex linear, i.e., the Christoffel form is complex bilinear.  
\end{prop}

\proof
First assume that the metric is K\"ahler. Then for any 3 vector fields
$X$, $Y$ and $Z$ we have
$$d\omega(X,Y,Z)=X\omega(Y,Z)-\omega([X,Y],Z)-\omega(Y,[X,Z])=0.$$
Here $[X,Y]$ denotes the commutator of the vector fields $X$ and $Y$.
Fix the values of $X$, $Y$ and $Z$ at some point $p$. We can always arrange
that $\nabla_Y X=\nabla_Z X=0$ at $p$ by changing $X$ in a neighborhood
of $p$. Then $[X,Y]=\nabla_X Y-\nabla_Y X=\nabla_X Y$. Analogously,
$[X,Z]=\nabla_X Z$. Recall that $\nabla_X$ depends on $X(p)$ only
(not on the derivatives of $X$). Finally, we obtain
$$X\omega(Y,Z)=\omega(\nabla_X Y,Z)+\omega(Y,\nabla_X Z).$$
On the other hand, by the compatibility of $\nabla$ with $g$,
$$X\omega(Y,Z)=Xg(iY,Z)=g(\nabla_X (iY),Z)+g(iY,\nabla_X Z).$$
Comparing our equations, we conclude that $g(\nabla_X(iY),Z)=
g(i\nabla_X Y,Z)$. Since $Z$ is arbitrary, $\nabla_X(iY)=i\nabla_X Y$,
i.e., the covariant differentiation is complex linear.

The above argument can be reversed. If the covariant differentiation
is complex linear, then $d\omega$ vanishes for all $X$, $Y$ and $Z$
such that $\nabla_Y X=\nabla_Z X=0$ at some given point $p$.
Since $X(p)$, $Y(p)$ and $Z(p)$ can take arbitrary values,
$d\omega=0$ at $p$. But $p$ is also arbitrary. Hence $d\omega=0$
everywhere. $\Box$

{\sc Proof of Proposition \ref{exp-hol}.}
Indeed, the second differential of an exponential map coincides
with $\Gamma$, but the latter is complex bilinear by Proposition
\ref{K-cri}. $\Box$

\section*{Appendix 2: Fubini spaces}

Fubini spaces are complex analogs of the classical geometries (Euclidean,
Riemann, Lobachevsky).

Consider the complex space $\C^{n+1}$ equipped with the pseudo-Hermitian form
$$H=Z_0\bar Z_0+\alpha\sum_{j=1}^n Z_j\bar Z_j$$
where $\alpha$ is some real number. 
The {\em pseudosphere} is a hypersurface $S$ given by the equation $H=1$.
Note that the pseudosphere is stable under the multiplication by
complex numbers with absolute value 1, i.e., under the scalar $U(1)$-action.
The quotient space $F=S/U(1)$ is called a {\em Fubini space}.

Denote by $C$ the cone where $H>0$. Then the Fubini space can be also
defined as the quotient $C/\C^*$
(since the intersection of a $\C^*$-orbit with $S$ is exactly a $U(1)$-orbit).
Hence, for $\alpha>0$ we obtain the
complex projective space $\CP^n$, for $\alpha=0$ --- the affine space $\C^n$
and for $\alpha<0$ --- the complex hyperbolic space $\H^n$.

Let us introduce a Riemannian metric in a Fubini space.
Suppose first that $\alpha\ne 0$.
Then $H$ induces a metric on $S$ (for $\alpha<0$, this metric will be
negative so we should take it with sign minus) which is stable under the
$U(1)$-action. Hence a Fubini space $F$ also inherits some metric.
Namely, the distance between $U(1)$-orbits is defined as the minimal
distance from a point of one orbit to a point of the other orbit. 
For $\alpha=0$, we should take the standard Euclidean metric on $\C^n=F$.

To get an affine model of a Fubini space $F$, it is enough to project
it to the hyperplane $\{Z_0=1\}$. Namely, each point $x\in F$ can be viewed
as a complex line in $C$. Take the intersection of this line with
$\{Z_0=1\}$. Under this projection, $F$ gets mapped to the whole hyperplane
(for $\alpha\geq 0$) or to the interior of a ball (for $\alpha<0$).
In particular, for $\alpha>0$ we get an affine chart of $\CP^n$. Metrics
of Fubini spaces written down in the affine models are called the
{\em Fubini metrics} on (parts of) $\C^n$. 

Let us deduce the coordinate expressions of Fubini metrics for $\alpha>0$.
Take a vector $v\in T_xF$ at some point $x\in F$. Consider a lift
$X$ of $x$ to $C$ and a lift $V$ of $v$ looking out of $X$. We can always
assume that $|X|^2=1$, i.e., $X\in S$
(all norms and inner products are with respect
to the form $H$). Denote by $W$ the projection of $V$ to the orthogonal
complement of $X$. Then the length of $v$ with respect to the Fubini metric
equals to the length of $W$ with respect to $H$:
$$|v|^2=\<W,W\>=\<V-\<V,X\>X,V-\<V,X\>X\>=
\<V,V\>-\<X,V\>\<V,X\>.$$
Now if $X$ is arbitrary (not necessarily of unit length), then the
formula for $|v|^2$ can be recovered by the homogeneity:
$$|v|^2=\frac{\<V,V\>\<X,X\>-\<X,V\>\<V,X\>}{\<X,X\>}.$$
The vector $X$ can be regarded as the collection of homogeneous coordinates
of the point $x$. In order to pass to affine coordinates, it is enough
to put $X_0=1$, $V_0=0$ ($X_0$ and $V_0$ stand for zero-coordinates of
$X$ and $V$ respectively). For $\alpha<0$ the above expression is to be
taken with negative sign.

\begin{prop}
\label{Fcirc}
All complex lines in Fubini metrics are totally geodesic surfaces. All
geodesics are (parts of) circles.
\end{prop}

\proof Note that a Fubini metric is preserved under the action of
(rather large) group of all $H$-unitary projective transformations.
Each complex line is stable under a one-parametric subgroup of rotations
around it. It follows that each complex line is a geodesic submanifold.
On a coordinate line passing through the origin we have a classical
geometry (standard Euclidean if $\alpha=0$, spherical in central projection
if $\alpha>0$ or Lobachevsky in the Poincar\'e disk model if $\alpha<0$).
Clearly all geodesics inside this line are circles. Any complex line can
be mapped to any other by an isometry. This concludes the proof. $\Box$

A useful characterization of Fubini spaces was given by Bochner \cite{Bo}.

First recall the definition of holomorphic sectional curvature. Let
$g$ be a Hermitian metric in an open subset $\Omega$ of $\C^n$.
Take a point $p\in\Omega$ and a vector $\xi$ going out from this point.
The vector $\xi$ defines a germ of complex line. Consider the
image of this germ under the exponential map of $g$.
The image is a germ of 2-dimensional surface at the point $p$.
Its Gauss curvature at $p$ is denoted by $K(\xi)$ and is called
the {\em holomorphic sectional curvature}. A metric is said to have
{\em constant holomorphic sectional curvature} if $K(\xi)$
depends neither on the direction of $\xi$ nor on the point $p$.

\begin{thm}[Bochner]
\label{curv}
A K\"ahler metric $g$ has constant holomorphic sectional curvature if
and only if $g$ is locally equivalent to a Fubini space via a holomorphic
change of variables.
\end{thm}

\section*{Appendix 3: Complexified Beltrami's theorem}

Proposition \ref{B} is a complexified version of classical Beltrami's
theorem \cite{Bel}: if all geodesics are parts of straight lines, then the
metric
is locally equivalent to Euclidean, Riemann or Lobachevsky. This complex
version can be deduced from the results of Bochner \cite{Bo}, Otsuki
and Tashiro \cite{OT}. Here we recall these results and also sketch another
proof of Proposition \ref{B} in dimension 2 which does not involve
curvature considerations.

\begin{lemma}
\label{pro}
Let $\Gamma:\C^n\to\C^n$ be a homogeneous polynomial of degree 2
over reals (i.e., not necessarily holomorphic). Suppose that $\Gamma(v)$
is everywhere proportional to $v$ with some complex coefficient $L(v)$.
Then $L$ is a complex-valued $\R$-linear function.
\end{lemma}

\proof
The coefficient $L$ is a complex-valued function defined everywhere except
perhaps 0. Since $\Gamma$ is quadratic over reals, it satisfies
the relation
$$\Gamma(v+w)+\Gamma(v-w)=2(\Gamma(v)+\Gamma(w))$$
for all $v,w\in\C^n$. Substituting $L(u)u$ for $\Gamma(u)$, we obtain:
$$v(L(v+w)+L(v-w)-2L(v))+w(L(v+w)-L(v-w)-2L(w))=0.$$
We can choose $v$ and $w$ to be linearly independent, so
$$L(v+w)+L(v-w)=2L(v),\quad L(v+w)-L(v-w)=2L(w).$$
If $v$ and $w$ are linearly dependent, this is also true due to
the homogeneity of $L$. Hence the equations above hold for all
$v$ and $w$. They imply that $L$ is $\R$-linear. $\Box$

\begin{prop}
\label{Gamma}
Consider a Hermitian metric $g$ in an open subset $\Omega$ of $\C^n$.
If all germs of complex lines lying in $\Omega$ are totally geodesic
submanifolds, then the Christoffel form is equal to $\Gamma(v)=L(v)v$ where
$L$ is some complex linear functional (for a definition of the
Christoffel form see Appendix 1).
\end{prop}

\proof
Consider an arbitrary vector $v$ at some point $x\in\Omega$
and a geodesic $\gamma$ passing through $x$ with velocity $v$.
By the equation of geodesics, $\ddot\gamma+\Gamma(v,v)=0$. But
since the geodesic lies in some complex line, $\ddot\gamma$ is
proportional to $\dot\gamma=v$ with some complex coefficient.
Therefore, $\Gamma(v,v)$ is proportional to $v$. By lemma \ref{pro},
$\Gamma(v,v)=L(v)v$. $\Box$

\begin{cor}
\label{Kah}
Under the assumptions of Proposition \ref{Gamma} the metric $g$ is K\"ahler. 
\end{cor}

\proof
By lemma \ref{Gamma}, $\Gamma(v,v)=L(v)v$. By the symmetry of
$\Gamma$, we have $\Gamma(v,w)=\frac 12(L(v)w+L(w)v)$. 

Now we can use the Hermitian property: $g(X,X)=g(iX,iX)$.
Apply the covariant differentiation $\nabla_Y$ to both sides
of this relation:
$$g(\nabla_YX,X)=g(\nabla_Y(iX),iX)=-g(i\nabla_Y(iX),X).$$
The standard connection $\nabla^0$ is complex linear, hence
$g(i\Gamma(Y,iX)+\Gamma(Y,X),X)=0$. Since $Y$ is arbitrary, it follows that
$L(iX)=iL(X)$, i.e., $L$ is complex linear. This means that
$\Gamma$ is $\C$-bilinear. By Proposition \ref{K-cri}, $g$ is
K\"ahler in this case. $\Box$

Consider 2 Hermitian
metrics $g'$ and $g''$ and denote the corresponding Levi-Civita connections
by $\nabla'$ and $\nabla''$ respectively. Recall that the difference
$\Gamma(X)=\nabla''_XX-\nabla'_XX$ is a vector-valued quadratic form. The
metrics
$g'$ and $g''$ are called {\em holomorphically projectively equivalent}
if $\Gamma(x)=L(x)x$ for all vectors $x$ where $L$ is a complex linear
functional (depending on point). Lemma \ref{Gamma} shows that if all
complex lines are geodesic surfaces, then the metric is holomorphically
projectively equivalent to the standard (flat) metric. 

Otsuki and Tashiro proved \cite{OT}
that a Hermitian metric that is holomorphically projectively equivalent
to a Fubini metric, has constant holomorphic sectional curvature. By
Bochner's theorem it is isometric to a Fubini space, an isometry being
a holomorphic map. But a holomorphic map taking (locally) complex lines
to complex lines is a complex projective transformation. Thus we obtain
Proposition \ref{B}. Below we sketch a more straight-forward proof of it
in dimension 2.

Suppose a metric $g$ in an open subset $\Omega$ of $\C^2$ satisfies the
conditions of Proposition \ref{B}. Choose a pair of constant linearly
independent vector fields $X$ and $Y$ in $\Omega$ and compose the Gram
determinant
$$G=G(X,Y)=\<X,X\>\<Y,Y\>-\<X,Y\>\<Y,X\>$$
where $\<,\>$ is the Hermitian inner product corresponding to $g$. 
Note that $G$ does not essentially depend on $X$ and $Y$. In fact,
it is well defined as a function up to a positive constant factor. 

\begin{lemma}
The Hermitian metric $h=g/G^{2/3}$ is constant along each complex line.
This means that for any vectors $v$ and $w$ of the same Euclidean length
lying in the same complex line we have $h(v)=h(w)$. 
\end{lemma}

\proof Let $X$ and $Y$ be constant linearly independent vector fields in
$\Omega$. Then we have $XG(X,Y)=3\Re L(X)G(X,Y)$. On the other hand,
$Xg(X,X)=2\Re L(X)g(X,X)$. It follows that $Xh=0$. Since $X$ is an arbitrary
constant vector field, $h$ must be constant along any real line. It remains
to note that a Hermitian metric constant along any real line is also
constant along any complex line. $\Box$

Note that $g$ can be recovered by $h$. Namely, if $H$ is the Gram determinant
of $h$, then $g=h/H^2$. It remains to describe all Hermitian metrics
that are constant along any complex line. This is not difficult to
accomplish. Any such metric considered as a function of a point $x$
and a vector $v$ out of $x$ is a second degree polynomial in
``complex momentum'' $v$ and ``complex angular momentum'' $x\wedge v$
(the wedge product is over complex numbers). One can readily verify that
these metrics provide Fubini metrics (modulo complex projective
transformations) after division by the square of their Gram determinants.


\begin{thebibliography}{9}

\bibitem{Kh} Khovanskii A.G., {\em Rectification of circles}, Sib. Mat. Zh.,
{\bf 21} (1980), 221--226
 
\bibitem{Iz} Izadi F.A., {\em Rectification of circles, spheres, and 
classical geometries}, PhD thesis, University of Toronto, (2001)

\bibitem{Tim}
Timorin V.A. {\em ``Rectification of circles and quaternions''},
Preprint\\ 
{\sf http://xxx.lanl.gov/math.DG/abs/0110144}.

\bibitem{Bo}
Bochner S. {\em ``Curvature in Hermitian metric''}, Bull. Amer. 
Math. Soc., {\bf 53} (1947), 149--195

\bibitem{OT}
Otsuki T., Tashiro Y. {\em ``On Curves in Kaehlerian spaces''},
Math. J. Okayama univ., {\bf 4} (1954), 57--78

\bibitem{Bel}
Beltrami E. {\em ``Risoluzione del problema: Riportare i punti di una
superficie sopra un piano in modo che le linee geodetiche vengano
rappresentate da linee rette''} Ann. Mat. pura appl., ser. 1, {\bf 7},
185--204

\end{thebibliography}
\end{document}